\documentclass[11pt,twocolumn,doublespace]{article}
\usepackage{graphicx,pslatex,amsmath,amsfonts,amssymb,amsthm}

\theoremstyle{remark}
\newtheorem{example}{Example}

\begin{document}


\parindent 0mm
\parskip 6pt


\title{Random Dynamical Systems}

\author{V\'\i tor Ara\'ujo\\Centro de Matem\'atica da
  Universidade do Porto\\Rua do Campo Alegre, 687, 4169-007
  Porto, Portugal\\E-mail: vdaraujo@fc.up.pt}


\maketitle


\section{Introduction}

The concept of random dynamical system is a comparatively recent
development combining ideas and methods from the well developed areas
of probability theory and dynamical systems.

Let us consider a mathematical model of some physical
process given by the iterates
$T_0^k=T_0\circ\stackrel{k}{\cdots}\circ T_0, \, k\ge1$, of a
smooth transformation $T_0:M\circlearrowleft$ of a manifold
into itself. A realization of the process with initial
condition $x_0$ is modelled by the sequence
$(T_0^k(x_0))_{k\ge1}$, the \emph{orbit} of $x_0$.

Due to our inaccurate knowledge of the particular physical
system or due to computational or theoretical limitations
(lack of sufficient computational power, inefficient
algorithms or insufficiently developed mathematical or
physical theory, for example), the mathematical models never
correspond exactly to the phenomenon they are meant to model.
Moreover  when considering practical systems we cannot avoid either
external noise or measurement or \emph{inaccuracy} errors,
so every realistic mathematical model should allow for small
errors along orbits not to disturb too much the long term
behavior.  To be able to cope with unavoidable uncertainty
about the ``correct'' parameter values, observed initial
states and even the specific mathematical formulation
involved, we let randomness be embedded within the model to
begin with.

We present the most basic classes of models in what follows, then
define the general concept and present some developments and
examples of applications.

\section{Dynamics with noise}
\label{sec:dynamics-with-noise}

To model random perturbations of a transformation $T_0$ we
may consider a transition from the image $T_0(x)$ to some
point according to a given probability law, obtaining a
Markov Chain, or, if $T_0$ depends on a parameter $p$, we
may choose $p$ at random at each iteration, which also can
be seen as a Markov Chain but whose transitions are strongly
correlated.

\subsection{Random noise}
\label{sec:random-noise}

Given $T_0:M\circlearrowleft$ and a family $\{p(\cdot\mid
x): x\in M\}$ of probability measures on $M$ such that the
support of $p(\cdot\mid x)$ is close to $T_0(x)$, the
\emph{random orbits} are sequences $(x_k)_{k\ge1}$ where
each $x_{k+1}$ is a random variable with law $p(\cdot\mid
x_k)$. This is a Markov Chain with state space $M$ and
transition probabilities $\{p(\cdot\mid x)\}_{ x\in M}$.
To extend the concept of invariant measure of a
transformation to this setting, we say that a probability
measure $\mu$ is \emph{stationary} if $\mu(A)=\int p(A\mid
x)\, d\mu(x)$ for every measurable (Borel) subset $A$. This
can be conveniently translated by saying that the
skew-product measure $\mu\times p^{\mathbb N}$ on $M\times
M^{\mathbb N}$ given by
\begin{eqnarray*}
&&
d(\mu\times p^{\mathbb N})(x_0,x_1,\dots, x_n,\dots)
\\
&&
=d\mu(x_0) p(dx_1\mid x_0)\cdots p(dx_{n+1}\mid x_n)\cdots
\end{eqnarray*}
is invariant by the shift map $\mathcal S:M\times M^{\mathbb
  N} \circlearrowleft$ on the space of orbits. Hence we may
use the Ergodic Theorem and get that time averages of every
continuous observable $\varphi:M\to\mathbb R$, i.e. writing
$\underline x=(x_k)_{k\ge0}$ and
\begin{eqnarray*}
  \tilde\varphi(\underline x) &=&
\lim_{n\to+\infty}\frac1n\sum_{k=0}^{n-1} \varphi(x_k)
\\
&=&
\lim_{n\to+\infty}\frac1n\sum_{k=0}^{n-1}
\varphi(\pi_0(\mathcal S^k (\underline x)))
\end{eqnarray*}
exist for $\mu\times p^{\mathbb N}$-almost all sequences
$\underline x$, where $\pi_0:M\times M^{\mathbb N}\to M$ is
the natural projection on the first coordinate. It is well
known that stationary measures always exist if the
transition probabilities $p(\cdot\mid x)$ depend
continuously on $x$.

A function $\varphi:M\to\mathbb R$ is \emph{invariant} if
$\varphi(x)=\int \varphi(z) p(dz\mid x)$ for $\mu$-almost
every $x$.  We then say that $\mu$ is \emph{ergodic} if
every invariant function is constant $\mu$-almost
everywhere. Using the Ergodic Theorem again, if $\mu$ is
ergodic, then $\tilde\varphi=\int\varphi\, d\mu$,
$\mu$-almost everywhere.

Stationary measures are the building blocks for more
sophisticated analysis involving e.g. asymptotic sojourn
times, Lyapunov exponents, decay of correlations, entropy
and/or dimensions, exit/entrance times from/to subsets of
$M$, to name just a few frequent notions of dynamical and
probabilistic/statistical nature.

\begin{example}[Random jumps]
\label{ex:random-jumps}
Given $\epsilon>0$ and $T_0:M\to M$, let us define
\[
p^\epsilon(A\mid x)=\frac{m(A\cap
  B(T_0(x),\epsilon))}{m(B(T_0(x),\epsilon))}
\]
where $m$ denotes some choice of Riemannian volume form on
$M$. Then $p^\epsilon(\cdot\mid x)$ is the normalized volume
restricted to the $\epsilon$-neighborhood of $T_0(x)$. This
defines a family of transition probabilities allowing the
points to ``jump'' from $T_0(x)$ to any point in the
$\epsilon$-neighborhood of $T_0(x)$ following a uniform
distribution law.
\end{example}

\subsection{Random maps}
\label{sec:random-maps}

Alternatively we may choose maps $T_1, T_2, \dots, T_k$
independently at random near $T_0$ according to a
probability law $\nu$ on the space $T(M)$ of maps, whose
support 
is close to $T_0$ in some
topology, and consider sequences $x_k=T_k\circ\cdots\circ
T_1 (x_0)$ obtained through random iteration, $k\ge1, \,
x_0\in M$.

This is again a Markov Chain whose transition probabilities
are given for any $x\in M$ by
\[
p(A\mid x)=\nu\big(\{ T\in T(M): T(x)\in A\}\big),
\]
so this model may be reduced to the first one. However in
the random maps setting we may associate to each random
orbit a sequence of maps which are iterated, enabling us to
use \emph{robust  properties} of the transformation
$T_0$ (i.e. properties which are known to hold for $T_0$ and
for every nearby map $T$) to derive properties of the
random orbits.

Under some regularity conditions on the map $x\mapsto
p(A\mid x)$ for every Borel subset $A$, it is possible to
represent random noise by random maps on suitably chosen
spaces of transformations. In fact the transition
probability measures obtained in the random maps setting
exhibit strong spatial correlation: $p(\cdot\mid x)$ is
close to $p(\cdot\mid y)$ is $x$ is near $y$.

If we have a parameterized family $T:\mathcal U\times M\to
M$ of maps we can specify the law $\nu$ by giving a
probability $\theta$ on $\mathcal U$. Then to every sequence
$T_1,\dots, T_k, \dots$ of maps of the given family we
associate a sequence $\omega_1,\dots, \omega_k,\dots$ of
parameters in $\mathcal U$ since
\[
T_k\circ\cdots\circ T_1=T_{\omega_k}\circ\dots\circ
T_{\omega_1}= T^k_{\omega_1,\dots,\omega_k}
\]
for all $k\ge1$, where we write $T_\omega(x)=T(\omega,x)$.
In this setting the shift map $\mathcal S$ becomes a
skew-product transformation
\[
\mathcal S: M\times \mathcal U^{\mathbb N}
\circlearrowleft\quad
(x, \underline \omega)\mapsto
\big(T_{\omega_1}(x),\sigma(\underline \omega)\big),
\]
to which many of the standard methods of dynamical
systems and ergodic theory can be applied, yielding stronger
results that can be interpreted in random terms.

\begin{example}[Parametric noise]
\label{ex:parametrical-noise}
Let $T: P\times M\to M$ be a smooth map where $P,M$ are
finite dimensional Riemannian manifolds. We fix $p_0\in P$,
denote by $m$ some choice of Riemannian volume form on $P$,
set $T_w(x)=T(w,x)$ and for every $\epsilon>0$ write
$\theta_\epsilon=(m(B(p_0,\epsilon))^{-1}\cdot(m\mid
B(p_0,\epsilon))$, the normalized restriction of $m$ to the
$\epsilon$-neighborhood of $p_0$. Then $(T_w)_{w\in P}$
together with $\theta_\epsilon$ defines a random
perturbation of $T_{p_0}$, for every small enough
$\epsilon>0$.
\end{example}

\begin{example}[Global additive perturbations]
\label{ex:global-additive}
  Let $M$ be a homogeneous space, i.e., a compact connected
   Lie Group admitting an invariant
  Riemannian metric. Fixing a neighborhood $\mathcal U$ of
  the identity $e\in M$ we can define a map $T:\mathcal U\times
  M\to M, (u,x)\mapsto L_u( T_0(x))$, where $L_u(x)=u\cdot
  x$ is the left translation associated to $u\in M$. The
  invariance of the metric means that left (an also right)
  translations are isometries, hence fixing $u\in \mathcal U$ and
  taking any $(x,v)\in TM$ we get
\begin{eqnarray*}
\|DT_u(x)\cdot v\|&=&\|DL_u(T_0(x))
(DT_0(x)\cdot v)\|
\\
&=&\|DT_0(x)\cdot v\|.
\end{eqnarray*}
In the particular case of $M={\mathbb T}^d$, the $d$-dimensional
torus, we have $T_u(x)=T_0(x)+u$ and this simplest case
suggests the name \emph{additive random perturbations} for
random perturbations defined using families of maps of this
type. 

For the probability measure on $\mathcal U$ we may take
$\theta_\epsilon$ any probability measure supported in the
$\epsilon$-neighborhood of $e$ and
absolutely continuous with respect to the Riemannian metric
on $M$, for any $\epsilon>0$ small enough.
\end{example}

\begin{example}[Local additive perturbations]
\label{ex:local-additive}
  If $M=\mathbb R^d$ and $U_0$ is a bounded open subset of
  $M$ strictly invariant under a diffeomorphism $T_0$,
  i.e., ${\rm closure\,}(T_0(U_0))\subset U_0$, then we
  can define an isometric random perturbation setting
  \begin{itemize}
  \item $V=T_0(U_0)$ (so that ${\rm closure\,}(V)={\rm
      closure\,}(T_0(U_0))\subset U_0$);
  \item $G\simeq{\mathbb R}^d$ the group of translations of
  $\mathbb R^d$;
\item $\mathcal V$ a small enough neighborhood of $0$
  in $G$.
  \end{itemize}
  Then for $v\in\mathcal V$ and $x\in V$ we set
  $T_v(x)=x+v$, with the standard notation for
  vector addition, and clearly $T_v$ is an isometry.
  For $\theta_\epsilon$ we may take any probability measure
  on the $\epsilon$-neighborhood of $0$, supported in
  $\mathcal V$ and absolutely continuous with respect to the
  volume in $\mathbb R^d$, for every small enough
  $\epsilon>0$.
\end{example}

\subsection{Random perturbations of flows}
\label{sec:perturbations-flows}

In the continuous time case the basic model to start with is
an ordinary differential equation $dX_t= f(t, X_t)dt$, where
$f:[0,+\infty)\to\mathcal X(M)$ and $\mathcal X(M)$ is the
family of vector fields in $M$. We embed randomness in the
differential equation basically through \emph{diffusion},
the perturbation is given by \emph{white noise} or
\emph{Brownian motion} ``added'' to the ordinary solution.

In this setting, assuming for simplicity that $M=\mathbb
R^n$, the random orbits are solutions of stochastic
differential equations
\[
d X_t = f(t, X_t) dt + \epsilon\cdot\sigma(t,X_t) dW_t,
\,\,\, 0\le t \le T, \, X_0=Z,
\]
where $Z$ is a random variable, $\epsilon,T>0$ and both
$f:[0,T]\times\mathbb R^n \to \mathbb R^n$ and
$\sigma:[0,T]\times\mathbb R^n \to \mathcal L(\mathbb
R^k,\mathbb R^n)$ are measurable functions. We have written
$\mathcal L(\mathbb R^k,\mathbb R^n)$ for the space of
linear maps $\mathbb R^k \to \mathbb R^n$ and $W_t$ for the
\emph{white noise process} on $\mathbb R^k$. The solution of
this equation is a stochastic process
\[
X:\mathbb R\times\Omega \to M,
\quad
(t,\omega)\mapsto X_t(\omega),
\]
for some (abstract) probability space $\Omega$, given by
\[
X_t=Z+\int_0^T f(s,X_s)\, ds +
\int_0^T \epsilon\cdot \sigma(s,X_s) \,dW_s,
\]
where the last term is a stochastic integral in the sense of
It\^o. Under reasonable
conditions on $f$ and $\sigma$, there exists a unique
solution \emph{with continuous paths}, i.e.
\[
[0,+\infty)\ni t \mapsto X_t(\omega)
\]
is continuous for almost all $\omega\in\Omega$ (in general
these paths are \emph{nowhere differentiable}).

Setting $Z=\delta_{x_0}$ the probability measure concentrated
on the point $x_0$, the initial point of the path is $x_0$
with probability $1$. We write $X_t(\omega)x_0$ for paths of
this type. Hence $x\mapsto X_t(\omega)x$ defines a map
$X_t(\omega):M\circlearrowleft$ which can be shown to be a
homeomorphism and even a diffeomorphisms under suitable
conditions on $f$ and $\sigma$. These maps satisfy a cocycle
property
\begin{eqnarray*}
  X_0(\omega) &=& {\mathrm Id}_M\,\,(\mbox{identity map
  of} \,\,  M),
\\
X_{t+s}(\omega) &=&
X_t(\theta(s)(\omega))\circ X_s(\omega),
\end{eqnarray*}
for $s,t\geq0$ and $\omega\in\Omega$, for a family of
measure preserving transformations
$\theta(s):(\Omega,\mathbb P)\circlearrowleft$ on a suitably
chosen probability space $(\Omega,\mathbb P)$. This enables
us to write the solution of this kind of equations also as a
skew-product.

\subsection{The abstract framework}
\label{sec:abstract-framework}

The illustrative particular cases presented can all be
written in skew-product form as follows. 

Let $(\Omega,\mathbb P)$ be a given probability space, which
will be the model for the noise, and $\mathbb T$ be time,
which usually means $\mathbb Z_+, \mathbb Z$ (discrete,
resp. invertible system) or $\mathbb R_+,\mathbb R$
(continuous, resp. invertible system).
A random dynamical system is a skew-product
\[
\mathcal S_t: \Omega\times M \circlearrowleft,
\,\,
(\omega, x)\mapsto (\theta(t)(\omega), \varphi(t,\omega)(x)),
\]
for all $t\in\mathbb T$, where $\theta:\mathbb
T\times\Omega\to \Omega$ is a family of measure preserving
maps $\theta(t):(\Omega,\mathbb P) \circlearrowleft$ and
$\varphi:\mathbb T\times\Omega\times M\to M$ is a family of
maps $\varphi(t,\omega): M\circlearrowleft$ satisfying the
cocycle property: for $s,t\in\mathbb T$, $\omega\in\Omega$
\begin{eqnarray*}
  \varphi(0,\omega)
&=&
{\mathrm Id}_M,
\\
\varphi(t+s,\omega) 
&=&
\varphi(t,\theta(s)(\omega))\circ \varphi(s,\omega).
\end{eqnarray*}


In this general setting an invariant measure for the random
dynamical system is any probability measure $\mu$ on
$\Omega\times M$ which is $\mathcal S_t$-invariant for all
$t\in\mathbb T$ and whose \emph{marginal} is $\mathbb P$,
i.e. $\mu(\mathcal S_t^{-1}(U))=\mu(U)$ and
$\mu(\pi_\Omega^{-1}(U))=\mathbb P(U)$ for every measurable
$U\subset \Omega\times M$, respectively, with
$\pi_\Omega:\Omega\times M\to\Omega$ the natural projection.

\begin{example}
  \label{ex:measures}
  In the setting of the previous examples of random
  perturbations of maps, the product measure
  $\eta=\mathbb P\times\mu$ on $\Omega\times M$, with
  $\Omega=\mathcal U^{\mathbb N}$, $\mathbb
  P=\theta_\epsilon^{\mathbb N}$ and $\mu$ any stationary
  measure, is clearly invariant. However not all invariant
  measures are product measures of this type.
\end{example}

Naturally an invariant measure is ergodic if every $\mathcal
S_t$-invariant function is $\mu$-almost everywhere constant.
i.e. if $\psi:\Omega\times M\to\mathbb R$ satisfies
$\psi\circ\mathcal S_t=\psi\,\,\, \mu$-almost everywhere for
every $t\in\mathbb T$, then $\psi$ is $\mu$-almost
everywhere constant.


\section{Applications}
\label{sec:appl-exampl}

We avoid presenting well established applications of both probability
or stochastic differential equations (solution of boundary value
problems,
optimal stopping, stochastic control etc) and dynamical systems (all
sort of models of physical, economic or biological phenomena,
solutions of differential equations, control systems etc), focusing
instead on topics where the subject sheds new light on these areas.


\subsection{Products of random matrices and 
the Multiplicative Ergodic Theorem}
\label{sec:mult-ergod-theor}

The following celebrated result on products of random matrices has
far-reaching applications on dynamical systems theory.

Let $(X_n)_{n\ge0}$ be a sequence of independent and
identically distributed random variables on the probability
space $(\Omega,P)$ with values in $\mathcal L(\mathbb
R^k,\mathbb R^k)$ such that $E(\log^+\|X_1\|)<+\infty$,
where $\log^+ x=\max\{0,\log x\}$ and $\|\cdot\|$ is a given
norm on $\mathcal L(\mathbb R^k,\mathbb R^k)$.  Writing
$\varphi_n(\omega)=X_n(\omega)\circ\dots\circ X_1(\omega)$
for all $n\ge1$ and $\omega\in\Omega$ we obtain a cocycle.
If we set
\begin{eqnarray*}
&B&=\{(\omega,y)\in\Omega\times\mathbb R^k:
\lim_{n\to+\infty} \frac1n\log\|\varphi_n(\omega)y\|
\\
&&\mbox{exists and is finite or is}\,\,-\infty\},
\,\, \mbox{and}
\\
&\Omega^\prime&=\{ \omega\in\Omega:
(\omega,y)\in B\,\,\mbox{for all}\,\, y\in\mathbb R^k\},
\end{eqnarray*}
then $\Omega^\prime$ contains a subset
$\Omega^{\prime\prime}$ of full probability and there exist
random variables (which might take the value $-\infty$) $
\lambda_1\geq\lambda_2\geq\dots\geq\lambda_k $ with the
following properties.

(1) Let $I=\{ k+1=i_1>i_2>\dots>i_{l+1}=1\}$ be any $(l+1)$-tuple of
  integers and then we define
  \begin{eqnarray*}
\Omega_I=\{
\omega\in\Omega^{\prime\prime}:
\lambda_i(\omega)=\lambda_j(\omega), i_{h}>i,j\geq i_{h+1},\,\,
\mbox{and}
\\
\lambda_{i_{h}}(\omega)> \lambda_{i_{h+1}}(\omega)\,\,
\mbox{for all}\,\, 1<h<l
\}
  \end{eqnarray*}
the set of elements where the sequence $\lambda_i$ jumps exactly at
the indexes in $I$.
Then for $\omega\in\Omega_I$, $1<h\le l$
\[
\Sigma_{I,h}(\omega)=
\{
y\in\mathbb R^k:
\lim_{n\to+\infty}\frac1n
\log\|\varphi_n(\omega)\|\le
\lambda_{i_h}(\omega)
\}
\]
is a vector subspace with dimension $i_{h-1}-1$.

(2) Setting $\Sigma_{I,k+1}(\omega)=\{0\}$, then
\[
\lim_{n\to+\infty}\frac1n
\log\|\varphi_n(\omega)\|
= \lambda_{i_h}(\omega),
\]
for every $y\in \Sigma_{I,h}(\omega)\setminus
\Sigma_{I,h+1}(\omega)$.

(3) For all $\omega\in\Omega^{\prime\prime}$ there exists
    the matrix
\[
A(\omega)=
\lim_{n\to+\infty}
\Big[\big(\varphi_n(\omega)\big)^* \varphi_n(\omega)\Big]^{1/2n}
\]
whose eigenvalues form the set $\{e^{\lambda_i}: i=1,\dots,k\}$.

The values of $\lambda_i$ are the random Lyapunov characteristics and
the corresponding subspaces are analogous to random eigenspaces. If
the sequence $(X_n)_{n\ge0}$ is ergodic, then the Lyapunov
characteristics become non-random constants, but the Lyapunov
subspaces are still random.

We can easily deduce the Multiplicative Ergodic Theorem for measure
preserving differentiable maps $(T_0,\mu)$ on manifolds $M$ from this
result. We assume for simplicity that $M\subset\mathbb R^k$ and set
$p(A\mid x)=\delta_{T_0(x)}(A)=1$ if $T_0(x)\in A$ and $0$ otherwise.
Then the measure $\mu\times p^\mathbb N$ on $M\times M^\mathbb N$ is
$\sigma$-invariant (as defined in Section 2) and we have that
$\pi_0\circ\sigma=T_0\circ\pi_0$, where $\pi_0:M^\mathbb N\to M$ is
the projection on the first coordinate, and also $(\pi_0)_*(\mu\times
p^\mathbb N)=\mu$. Then setting for $n\ge1$
\[
\begin{array}[l]{llll}
X:&M&\to&\mathcal L(\mathbb R^k,\mathbb R^k)
\\
&  x & \mapsto & DT_0(x)
\end{array}
\mbox{and}\,\,\,
X_n=X\circ\pi_0\circ\sigma^n
\]
we obtain a stationary sequence to which we can apply the
previous result, obtaining the existence of Lyapunov
exponents and of Lyapunov subspaces on a \emph{full measure subset
for any $C^1$ measure preserving dynamical system}.

By a standard extension of the previous setup we obtain a
random version of the multiplicative ergodic theorem. We
take a family of skew-product maps $\mathcal
S_t:\Omega\times M\circlearrowleft$ as in
Subsection~\ref{sec:abstract-framework} with an invariant
probability measure $\mu$ and such that $\varphi(t,\omega):
M\circlearrowleft$ is (for simplicity) a local
diffeomorphism. We then consider the stationary family
\[
\begin{array}[l]{llll}
X_t:&\Omega&\to&\mathcal L(TM)
\\
& \omega & \mapsto & D\varphi(t,\omega): TM\circlearrowleft
\end{array},
\quad t\in\mathbb T,
\]
where $D\varphi(t,\omega)$ is the tangent map to
$\varphi(t,\omega)$. This is a cocycle since for all
$t,s\in\mathbb T,\, \omega\in\Omega$ we have
\[
X(s+t,\omega)=X(s,\theta(t)\omega)\circ X(t,\omega).
\]
If we assume that
\[
\sup_{0\le t\le 1} \, \sup_{x\in M} \, \big(\log^+ \|
D\varphi(t,\omega)(x) \|\big) \in L^1(\Omega,\mathbb P),
\]
where $\|\cdot\|$ denotes the norm on the corresponding
space of linear maps given by the induced norm (from the
Riemannian metric) on the appropriate tangent spaces, then
we obtain a sequence of random variables (which might take
the value $-\infty$) $
\lambda_1\geq\lambda_2\geq\dots\geq\lambda_k$, with $k$
being the dimension of $M$, such that
\[
\lim_{t\to+\infty}\frac1{t}\log\|X_t(\omega,x)y\|
=\lambda_{i}(\omega,x)
\]
for every $y\in
E_i\omega,x)=\Sigma_i(\omega,x)\setminus\Sigma_{i+1}(\omega,x)$
and $i=1,\dots,k+1$ where $(\Sigma_i(\omega,x))_i$ is a
sequence of vector subspaces in $T_{x}M$ as before,
measurable with respect to $(\omega,x)$.
In this setting the subspaces $E_i(\omega,x)$ and the
Lyapunov exponents are invariant, i.e. for all $t\in\mathbb
T$ and $\mu$ almost
every $(\omega,x)\in\Omega\times M$ we have
\[
\lambda_i(\mathcal S_t(\omega,x))=\lambda_i(\omega,x)\,\,\,
\mbox{and}\,\,\,
E_i(\mathcal S_t(\omega,x))=E_i(\omega,x).
\]

The \emph{dependence of Lyapunov exponents on the map} $T_0$
has been a fruitful and central research program in
dynamical systems for decades extending to the present day.
The random multiplicative ergodic theorem sets the stage for
the study of the stability of Lyapunov exponents under
random perturbations.


\subsection{Stochastic stability of physical measures}
\label{sec:stochastic-stability}

The development of the theory of dynamical systems has shown
that models involving expressions as simple as quadratic
polynomials (as the \emph{logistic family} or \emph{H\'enon
  attractor}), or autonomous ordinary differential equations
with a hyperbolic singularity of saddle-type, as the
\emph{Lorenz flow}, exhibit \emph{sensitive dependence on
  initial conditions}, a common feature of \emph{chaotic
  dynamics}: small initial differences are rapidly augmented
as time passes, causing two trajectories originally coming
from practically indistinguishable points to behave in a
completely different manner after a short while.  Long term
predictions based on such models are unfeasible since it is
not possible to both specify initial conditions with
arbitrary accuracy and numerically calculate with arbitrary
precision.

\subsubsection*{Physical measures}
\label{sec:physical-measures}

Inspired by an analogous situation of unpredictability faced
in the field of Statistical Mechanics/Thermodynamics,
researchers focused on the statistics of the data provided
by the time averages of some observable (a continuous
function on the manifold) of the system. Time averages are
guaranteed to exist for a positive volume subset of initial
states (also called an \emph{observable subset}) on the
mathematical model if the transformation, or the flow
associated to the ordinary differential equation, admits a
smooth invariant measure (a density) or a \emph{physical}
measure.

Indeed, if $\mu_0$ is an ergodic invariant measure for the
transformation $T_0$, then the Ergodic Theorem ensures that
for every $\mu$-integrable function $\varphi:M\to\mathbb R$
and for $\mu$-almost every point $x$ in the manifold $M$ the
time average $ \tilde\varphi(x)=\lim_{n\to+\infty}
n^{-1}\sum_{j=0}^{n-1} \varphi(T_0^j (x)) $ exists and
equals the space average $\int\varphi\,d\mu_0$. A
\emph{physical measure} $\mu$ is an invariant probability
measure for which it is \emph{required} that \emph{time
  averages of every continuous function $\varphi$ exist for
  a positive Lebesgue measure (volume) subset of the space
  and be equal to the space average $\mu(\varphi)$}.

We note that if $\mu$ is  a density, that is, is
absolutely continuous with respect to the volume measure,
then the Ergodic Theorem ensures that $\mu$ is physical.
However not every physical measure is absolutely continuous.
To see why in a simple example we just have to consider a
singularity $p$ of a vector field which is an attracting
fixed point (a sink), then the Dirac mass $\delta_p$
concentrated on $p$ is a physical probability measure, since
every orbit in the basin of attraction of $p$ will have
asymptotic time averages for any continuous observable
$\varphi$ given by $\varphi(p)=\delta_p(\varphi)$.

Physical measures need not be unique or even exist in
general, but when they do exist it is desirable that
\emph{the set of points whose asymptotic time averages are
  described by physical measures} (such set is called the
\emph{basin} of the physical measures) \emph{be of full
  Lebesgue measure} --- only an exceptional set of points
with zero volume would not have a well defined asymptotic
behavior. This is yet far from being proved for most
dynamical systems, in spite of much recent progress in this
direction. 

There are robust examples of systems admitting several
physical measures whose basins together are of full Lebesgue
measure, where \emph{robust} means that there are whole open
sets of maps of a manifold in the $C^2$ topology exhibiting
these features.    For typical parameterized families
of one-dimensional unimodal maps (maps of the circle or of
the interval with a unique critical point) it is known that
the above scenario holds true for Lebesgue almost every
parameter. It is known that there are systems
admitting no physical measure, but the only known cases are
not robust, i.e. there are systems arbitrarily close which
admit physical measures.

It is hoped that conclusions drawn from models admitting
physical measures to be effectively observable in the
physical processes being modelled. In order to lend more
weight to this expectation researchers demand stability
properties from such invariant measures.

\subsubsection*{Stochastic stability}
\label{sec:stochastic-stability-1}

There are two main issues when we are given a mathematical
model, both theoretical but with practical consequences. The
first one is to describe the asymptotic behavior of most
orbits, that is, to understand where do orbits go when time
tends to infinity. The second and equally important one is
to ascertain whether the asymptotic behavior is stable under
small changes of the system, i.e. whether the limiting
behavior is still essentially the same after small changes
to the evolution law. In fact since models are always
simplifications of the real system (we cannot ever take into
account the whole state of the universe in any model), the
lack of stability considerably weakens the conclusions drawn
from such models, because some properties might be specific
to it and not in any way resemblant of the real system.

Random dynamical systems come into play in this setting when
we need to check whether a given model is
stable under small random changes to the evolution law.

In more precise terms, we suppose that we are given a
dynamical system (a transformation or a flow) admitting a
physical measure $\mu_0$, and we take any random dynamical
system obtained from this one through the introduction of
small random perturbations on the dynamics, as in Examples
\ref{ex:random-jumps}-~\ref{ex:local-additive} or in
Subsection~\ref{sec:perturbations-flows}, with the noise
level $\epsilon>0$ close to zero.

In this setting if, for any choice $\mu_\epsilon$ of
invariant measure for the random dynamical system for all
$\epsilon>0$ small enough, the set of accumulation points of
the family $(\mu_\epsilon)_{\epsilon>0}$, when $\epsilon$
tends to $0$ --- also known as \emph{zero noise limits} --- is
formed by physical measures or, more generally, by  convex linear
combinations of physical measures, then the original
unperturbed dynamical system is \emph{stochastically
  stable}.

This intuitively means that the asymptotic behavior measured
through time averages of continuous observables for the
random system is close to the behavior of the unperturbed
system.

Recent progress in one-dimensional dynamics has shown that,
for typical families $(f_t)_{t\in(0,1)}$ of maps of the
circle or of the interval having a unique critical point, a
full Lebesgue measure subset $T$ of the set of parameters is
such that, for $t\in T$, the dynamics of $f_t$ admits a
unique stochastically stable (under additive noise type
random perturbations) physical measure $\mu_t$ whose basin
has full measure in the ambient space (either the circle or
the interval). Therefore models involving one-dimensional
unimodal maps typically are stochastically stable.

In many settings (e.g. low dimensional dynamical systems)
Lyapunov exponents can be given by time averages of
continuous functions --- for example the time average of
$\log\|DT_0\|$ gives the biggest exponent. In this case
stochastic stability directly implies stability of the
Lyapunov exponents under small random perturbations of the
dynamics.


\begin{example}[Stochastically stable examples]
\label{ex:stochstable-additive}
Let $T_0:\mathbb S^1\circlearrowleft$ be a map such that
$\lambda$, the Lebesgue (length) measure on the circle, is
$T_0$-invariant and ergodic. Then $\lambda$ is physical. 

We consider the parameterized family $T_t:\mathbb
S^1\times\mathbb S^1\to \mathbb S^1, (t,x)\mapsto x+t$
and a family of probability measures $\theta_\epsilon=
(\lambda(-\epsilon,\epsilon))^{-1}\cdot
(\lambda\mid(-\epsilon,\epsilon))$ given by the normalized
restriction of $\lambda$ to the $\epsilon$-neighborhood of
$0$, where we regard $\mathbb S^1$ as the Lie group
$\mathbb R/\mathbb Z$ and use additive notation for the
group operation. Since $\lambda$ is $T_t$-invariant for
every $t\in\mathbb S^1$, $\lambda$ is also an invariant
measure for the measure preserving random system
\[
\mathcal S: (\mathbb
S^1\times\Omega^{\mathbb N},\lambda\times\theta_\epsilon^{\mathbb
  N})\circlearrowleft,
\]
for every $\epsilon>0$, where $\Omega=(\mathbb S^1)^{\mathbb
  N}$. Hence $(T_0,\lambda)$ is stochastically stable under
additive noise perturbations.

Concrete examples can be irrational rotations,
$T_0(x)=x+\alpha$ with $\alpha\in\mathbb R\setminus\mathbb Q$, or
expanding maps of the circle, $T_0(x)=b\cdot x$ for some
$b\in\mathbb N,\, n\ge2$. Analogous examples exist in higher
dimensional tori.
\end{example}

\begin{example}[Stochastic stability depends on the type of noise]
  \label{ex:choicenoise}
  In spite of the straightforward way to obtain stochastic
  stability in Example~\ref{ex:stochstable-additive}, for
  e.g. an expanding circle map $T_0(x)=2\cdot x$, we can
  choose a continuous family of probability measures
  $\theta_\epsilon$ such that the same map $T_0$ is not
  stochastically stable. 
  
  It is well known that $\lambda$ is the unique
  absolutely continuous invariant measure for $T_0$ and also
  the unique physical measure.  Given $\epsilon>0$ small let
  us define transition probability measures as follows
\[
p_\epsilon(\cdot\mid z)=
\frac{\lambda\mid 
[\phi_\epsilon(z)-\epsilon,\phi_\epsilon(z)+\epsilon]}
{\lambda ([\phi_\epsilon(z)-\epsilon,\phi_\epsilon(z)+\epsilon])},
\]
where $\phi_\epsilon\mid(-\epsilon,\epsilon)\equiv 0$,
$\phi_\epsilon\mid [\mathbb
S^1\setminus(-2\epsilon,2\epsilon)]\equiv T_0$ and over
$(-2\epsilon,-\epsilon]\cup[\epsilon,2\epsilon)$ we define
$\phi_\epsilon$ by interpolation in order that it be smooth.

In this setting every random orbit starting at
$(-\epsilon,\epsilon)$ never leaves this neighborhood in the
future. Moreover it is easy to see that every random orbit
eventually enters $(-\epsilon,\epsilon)$. Hence every
invariant probability measure $\mu_\epsilon$ for this Markov
Chain model is supported in $[-\epsilon,\epsilon]$. Thus
letting $\epsilon\to0$ we see that the only zero-noise limit
is $\delta_0$ the Dirac mass concentrated at $0$, which is
not a physical measure for $T_0$.

This construction can be done in a random maps setting, but
only in the $C^0$ topology --- it is not possible to realize
this Markov Chain by random maps that are $C^1$ close to
$T_0$ for $\epsilon$ near $0$.
\end{example}


\subsection{Characterization of measures satisfying the Entropy
  Formula}
\label{sec:extend-known-prop}

A lot of work has been put in recent years in extending
important results from dynamical systems to the random
setting. Among many examples we mention the local conjugacy
between the dynamics near a hyperbolic fixed point and the
action of the derivative of the map on the tangent space,
the stable/unstable manifold theorems for hyperbolic
invariant sets and the notions and properties of metric and
topological entropy, dimensions and equilibrium states for
potentials on \emph{random (or fuzzy) sets}.

The characterization of measures satisfying the Entropy
Formula is one important result whose extension to the
setting of iteration of independent and identically
distributed random maps has recently had interesting new
consequences back into non-random dynamical systems.

\subsubsection*{Metric entropy for random perturbations}
\label{sec:metr-entr-rand}

Given a probability measure $\mu$ and a partition $\xi$ of
$M$, except perhaps for a subset of $\mu$-null measure, the
\emph{entropy of $\mu$ with respect to $\xi$} is defined to
be
\[
H_{\mu}(\xi)=-\sum_{R\in\xi} \mu(R)\log\mu(R)
\]
where we convention $0\log0=0$. Given another finite
partition $\zeta$ we write $\xi\vee\zeta$ to indicate the
partition obtained through intersection of every element of
$\xi$ with every element of $\zeta$, and analogously for any
finite number of partitions.  If $\mu$ is also a
stationary measure for a random maps model (as in
Subsection~\ref{sec:random-maps}), then for any finite
measurable partition $\xi$ of $M$,
$$
h_{\mu}(\xi) = \inf_{n\ge1} \frac{1}{n} \int H_{\mu}
\big( \bigvee_{i=0}^{n-1} (T^i_{\underline\omega})^{-1} (\xi) \big) d
p^{\mathbb N} (\underline\omega)
 $$
 is finite and is called \emph{the entropy of the random
 dynamical system} with respect to $\xi$ and to $\mu$.

We define $h_{\mu}= \sup_\xi \, h_{\mu}( \xi)$ as the
\emph{metric entropy} of the random dynamical system, where
the supremo is taken over all $\mu$-measurable partitions.
An important point here is the following notion: setting
$\mathcal A$ the Borel $\sigma$-algebra of $M$, we say that
a finite partition $\xi$ of $M$ is a \emph{random generating
  partition} for $\mathcal A$ if
\[
\bigvee_{i=0}^{+\infty} (T_{\underline\omega}^i)^{-1} (\xi)
=\mathcal A
\]
(except $\mu$-null sets) for $p^{\mathbb N}$-almost all
$\omega\in\Omega=\mathcal U^{\mathbb N}$. Then a classical
result from Ergodic Theory ensures that we can calculate the
entropy using only a random generating partition $\xi$: we
have $h_{\mu}=h_{\mu}( \xi)$.

\subsubsection*{The Entropy Formula}
\label{sec:entropy-formula}

There exists a general relation ensuring that the entropy of a
measure preserving differentiable transformation $(T_0,\mu)$
on a compact Riemannian manifold is bounded from above by
the sum of the positive Lyapunov exponents of $T_0$
\[
h_{\mu}(T_0)\le \int \sum_{\lambda_i(x)>0}
\!\lambda_i(x)\,\, d\mu(x).
\]
The equality (\emph{Entropy Formula}) was first shown to
hold for diffeomorphisms preserving a measure equivalent to
the Riemannian volume, and then the measures satisfying the
Entropy Formula were characterized: for \emph{$C^2$
  diffeomorphisms the equality holds if, and only if, the
  disintegration of $\mu$ along the unstable manifolds is
  formed by measures absolutely continuous with respect to
  the Riemannian volume restricted to those submanifolds}.
The \emph{unstable manifolds} are the submanifolds of $M$
everywhere tangent to the Lyapunov subspaces corresponding
to all positive Lyapunov exponents, the analogous to
``integrating the distribution of Lyapunov subspaces
corresponding to positive exponents'' --- this particular
point is a main subject of smooth ergodic theory for
non-uniformly hyperbolic dynamics.

Both the inequality and the characterization of stationary
measures satisfying the Entropy Formula were extended to
random iterations of independent and identically distributed
$C^2$ maps (non-injective and admitting critical points),
and the inequality reads
\[
h_{\mu}\le \int\!\!\int \sum_{\lambda_i(x,\omega)>0}\!
\lambda_i(x,\omega)\,\, d\mu(x)\,d p^{\mathbb N}(\omega).
\]
where the functions $\lambda_i$ are the random variables
provided by the Random Multiplicative Ergodic Theorem.


\subsection{Construction of physical measures
  as zero-noise limits}
\label{sec:constr-phys-meas}

The characterization of measures which satisfy the Entropy
Formula enables us to construct physical measures as
zero-noise limits of random invariant measures in some
settings, outlined in what follows, obtaining in the process
that the physical measures so constructed are also
stochastically stable.

The physical measures obtained in this manner arguably are
\emph{natural measures} for the system, since they are both
stable under (certain types of) random perturbations and
describe the asymptotic behavior of the system for a
positive volume subset of initial conditions.  This is a
significant contribution to the state-of-the-art of present
knowledge on Dynamics from the perspective of Random
Dynamical Systems.

\subsubsection*{Hyperbolic measures and the Entropy Formula}
\label{sec:hyperb-meas-entr}

The main idea is that an ergodic invariant measure $\mu$ for
a diffeomorphism $T_0$ which satisfies the Entropy Formula
and whose Lyapunov exponents are everywhere non-zero (known
as \emph{hyperbolic measure}) necessarily is a
\emph{physical measure} for $T_0$.  This follows from
standard arguments of smooth non-uniformly hyperbolic
ergodic theory.

Indeed $\mu$ satisfies the Entropy Formula if, and only if,
$\mu$ disintegrates into densities along the unstable
submanifolds of $T_0$. The unstable manifolds $W^u(x)$ are
tangent to the subspace corresponding to every positive
Lyapunov exponent at $\mu$-almost every point $x$, they are
an invariant family, i.e. $T_0(W^u(x))=W^u(x)$ for
$\mu$-almost every $x$, and distances on them are uniformly
contracted under iteration by $T_0^{-1}$.

If we know that the exponents along the complementary
directions are non-zero, then they must be negative and
smooth ergodic theory ensures that there exist \emph{stable
  manifolds}, which are submanifolds $W^s(x)$ of $M$
everywhere tangent to the subspace of negative Lyapunov
exponents at $\mu$-almost every point $x$, form a
$T_0$-invariant family ($T_0(W^s(x))=W^s(x)$, $\mu$-almost
everywhere), and distances on them are uniformly contracted
under iteration by $T_0$.

We still need to understand that time averages are
constant along both stable and unstable manifolds, and that
the families of stable and unstable manifolds are absolutely
continuous, in order to realize how an hyperbolic measure is
a physical measure.

Given $y\in W^s(x)$ the time averages of $x$ and $y$
coincide for continuous observables simply because ${\rm
  dist\,}(T_0^n(x),T_0^n(y))\to0$ when $n\to+\infty$. For
unstable manifolds the same holds when considering time
averages for $T_0^{-1}$. Since forward and backward time
averages are equal $\mu$-almost everywhere, we see that the
set of points having asymptotic time averages given by $\mu$
has positive Lebesgue measure if the following set
\[
B=\bigcup\{ W^s(y) : y\in W^u(x)\cap {\rm supp\,}(\mu) \}
\]
has positive volume in $M$, for some $x$ whose time averages
are well defined.

Now, stable and unstable manifolds are transverse everywhere
where they are defined, but they are only defined
$\mu$-almost everywhere and depend measurably on the base
point, so we cannot use  transversality arguments
from differential topology, in spite of $W^u(x)\cap {\rm
  supp\,}(\mu)$ \emph{having positive volume in $W^u(x)$} by
the existence of a smooth disintegration of $\mu$ along the
unstable manifolds. However it is known for smooth
($C^2$) transformations that the families of stable and unstable
manifolds are \emph{absolutely continuous}, meaning that
projections along leaves preserve sets of zero volume. This
is precisely what is needed for measure-theoretic arguments
to show that $B$ has positive volume.

\subsubsection*{Zero-noise limits satisfying the Entropy
  Formula}
\label{sec:zero-noise-limits}

Using the extension of the characterization of measures
satisfying the Entropy Formula for the random maps setting,
we can build random dynamical systems, which are small
random perturbations of a map $T_0$, having invariant
measures $\mu_\epsilon$ satisfying the Entropy Formula for
all sufficiently small $\epsilon>0$. Indeed it is enough to
construct small random perturbations of $T_0$ having
absolutely continuous invariant probability measures
$\mu_\epsilon$ for all small enough $\epsilon>0$.

In order to obtain such random dynamical systems we choose
families of maps $T:\mathcal U\times M\to M$ and of
probability measures $(\theta_\epsilon)_{\epsilon>0}$ as in
Examples~\ref{ex:global-additive}
and~\ref{ex:local-additive}, where we assume that
$o\in\mathcal U$ so that $T_0$ belongs to the family.
Letting $T_x(u)=T(u,x)$ for all $(u,x)\in\mathcal U\times
M$, we then have that $T_x(\theta_\epsilon)$ is absolutely
continuous.  This means that sets of perturbations of
positive $\theta_{\epsilon}$-measure send points of $M$ onto
positive volume subsets of $M$.  This kind of perturbation
can be constructed for every continuous map of any manifold.

In this setting we have that any invariant probability
measure for the associated skew-product map $\mathcal
S:\Omega\times M\circlearrowleft$ of the form
$\theta_\epsilon^{\mathbb N}\times\mu_\epsilon$ is such that
$\mu_\epsilon$ is absolutely continuous with respect to
volume on $M$. Then the Entropy Formula holds
\[
h_{\mu_\epsilon} = \int\!\!\int \sum_{\lambda_i(x,\omega)>0}\!
\lambda_i(x,\omega)\,\, d\mu_\epsilon(x)\, d
\theta_\epsilon^{\mathbb N}(\omega).
\]
Having this and knowing the characterization of measures
satisfying the Entropy Formula, it is natural to look for
conditions under which we can guarantee that the above
inequality extends to any zero-noise
limit $\mu_0$ of $\mu_\epsilon$ when $\epsilon\to0$. In that
case $\mu_0$ satisfies the Entropy Formula for $T_0$.

If in addition we are able to show that $\mu_0$ is a
hyperbolic measure, then we obtain a physical measure for
$T_0$ which is stochastically stable by construction.

These ideas can be carried out completely for hyperbolic
diffeomorphisms, i.e. maps admitting an continuous invariant
splitting of the tangent space into two sub-bundles $E\oplus
F$ defined everywhere with bounded angles, whose Lyapunov
exponents are negative along $E$ and positive along $F$.
Recently maps satisfying weaker conditions where shown to
admit stochastically stable physical measures following the
same ideas.

These ideas also have applications to the construction and
stochastic stability of physical measure for \emph{strange
  attractors} and for all mathematical models involving
ordinary differential equations or iterations of maps.



\section*{See also}
Equilibrium  statistical mechanics\\
Dynamical systems\\
Global analysis\\
Non-equilibrium statistical mechanics\\
Ordinary and partial differential equations\\
Stochastic methods\\
Strange attractors\\

\section*{Keywords}
Dynamical system\\
Flows\\
Orbits\\
Ordinary differential equations\\
Markov chains\\
Multiplicative Ergodic Theorem\\
Physical measures\\
Products of random matrices\\
Random maps\\
Random orbits\\
Random perturbations\\
Stochastic processes\\
Stochastic differential equations\\
Stochastic flows of diffeomorphisms\\
Stochastic stability\\


\section*{Further Reading}

L. Arnold, (1998),
\newblock {\em Random dynamical systems}.
\newblock Springer-Verlag, Berlin.

P.~Billingsley, (1965),
\newblock {\em Ergodic theory and information}.
\newblock J. Wiley \& Sons, New York.

P.~Billingsley, (1985),
\newblock {\em Probability and Measure}.
\newblock John Wiley and Sons, New York, 3rd edition.

J.~Doob, (1953),
\newblock {\em Stochastic Processes}.
\newblock Wiley, New York.

A.~Fathi, M.~Herman, and J.-C. Yoccoz (1983), \newblock A
proof of {P}esin's stable manifold theorem.  \newblock In
{\em Geometric dynamics (Rio de Janeiro 1981) edited by J.
  Palis}, volume 1007 of {\em Lect. Notes in Math.}, pages
177--215. Springer Verlag, New York.

Y. Kifer, (1986),
\newblock {\em Ergodic theory of random perturbations}.
\newblock Birkh{\"a}user, Boston.

Y. Kifer, (1988),
\newblock {\em Random perturbations of dynamical systems}.
\newblock Birkh{\"a}user, Boston.

H. Kunita, (1990),
\emph{Stochastic flows and stochastic differential
  equations},
Cambridge University Press, Cambridge.

F~Ledrappier and L.-S. Young, (1998).
\newblock Entropy formula for random transformations.
\newblock {\em Probab. Theory and Related Fields}, 80(2): 217--240.

B. {\O}skendal, (1992),
\newblock {\em Stochastic Differential Equations}.
\newblock Universitext. Springer-Verlag, Berlin, 3rd edition.

P.-D. Liu and M.~Qian, (1995)
\newblock {\em Smooth ergodic theory of random dynamical
  systems}, volume 1606 of {\em Lect. Notes in Math.}
\newblock Springer Verlag.

P.~Walters, (1982),
\newblock {\em An introduction to ergodic theory}.
\newblock Springer Verlag.

Bonatti C, Díaz L, Viana, M (2005)
\emph{Dynamics beyond uniform hyperbolicity. A global
geometric and probabilistic perspective.}
Encyclopaedia of Mathematical Sciences, 102. 
Mathematical Physics, III. 
Springer-Verlag, Berlin, 2005.

M. Viana (2000).
\newblock What's new on {L}orenz strange attractor.
\newblock {\em Mathematical Intelligencer}, 22(3): 6--19.



\end{document}